\documentclass{article}
\begin{document}
\newtheorem{proposition}{Proposition}[section]
\newtheorem{definition}{Definition}[section]
\newtheorem{lemma}{Lemma}[section]
 
\title{\bf Even (Odd) Roots and Simple Leibniz Algebras with Lie Factor $s\ell _2$}
\author{Keqin Liu\\Department of Mathematics\\The University of British Columbia\\Vancouver, BC\\
Canada, V6T 1Z2}
\date{December 27, 2006}
\maketitle

\begin{abstract}We introduce $3$-irreducible modules, even roots and odd roots for Leibniz algebras, produce a basis for a root space of a Leibniz algebra with a semisimple Lie factor, and classify finite dimensional  simple Leibniz algebras with Lie factor $s\ell _2$.
\end{abstract}

In this paper, all vector spaces are vector spaces over an algebraically closed field $\mathbf{k}$ of characteristic $0$. 

\medskip
This paper consists of three sections. In Section 1, we introduce $3$-irreducible modules which are natural building blocks among the modules over Leibniz algebras. In Section 2, we produce a basis of a weight space of a module over a right nilpotent Leibniz algebra by using the Extended Lie's Theorem, and construct the root space decomposition of a Leibniz algebra by using its Cartan subalgebra, even roots and odd roots. In Section 3, we classify finite dimensional  simple Leibniz algebras with Lie factor $s\ell _2$.

\section{$3$-Irreducible Modules}

We begin this section with the definition of a (right) Leibniz algebra (\cite{LT}).

\begin{definition}\label{def1.1} A vector space $L$ is called a {\bf (right) Leibniz algebra} if there exists a binary operation $\langle \, , \, \rangle$: $L\times L\to L$, called the {\bf angle bracket},  such that the {\bf (right) Leibniz identity} holds:
\begin{equation}\label{eq1.1}
\langle\langle x, y\rangle , z\rangle=\langle x, \langle y, z\rangle\rangle+ 
\langle\langle x,  z\rangle , y\rangle \qquad\mbox{for $x,y,z\in L$.}
\end{equation}
\end{definition}

A Leibniz algebra $L$ with an angle bracket $\langle \, , \, \rangle$ is also denoted by 
$(L, \, \langle \, , \, \rangle )$. 

\begin{definition}\label{def1.2} Let $I$ be a subspace of a Leibniz algebra 
$(L, \, \langle \, , \, \rangle )$.
\begin{description}
\item[(i)] $I$ is called a (Leibniz) {\bf subalgebra} of $L$ if 
$\langle I,  I\rangle\subseteq I$. 
\item[(ii)] $I$ is called an {\bf ideal} of $L$ if $\langle I, L\rangle\subseteq I$ and 
$\langle L, I\rangle\subseteq I$.
\end{description}
\end{definition}

The {\bf annihilator} $L^{ann}$ of a Leibniz algebra $L$ over a field $\bf{k}$ is defined by 
\begin{equation}\label{eq1.2}
L^{ann}:=\sum _{x\in L}\mathbf{k}\langle x, x\rangle.
\end{equation}
Clearly, $L^{ann}$ is an ideal of $L$ and
\begin{equation}\label{eq1.3}
\langle L ,  L^{ann} \rangle =0.
\end{equation}
$L^{ann}$ can also be expressed as follows:
\begin{equation}\label{eq1.4}
L^{ann}:=\sum _{x, y\in L}\mathbf{k}(\langle x, y\rangle +\langle y, x\rangle ).
\end{equation}
The quotient $L^{lie}:=\displaystyle\frac{L}{L^{ann}}$ is called the {\bf Lie-factor} of $L$. It is easy to check that $L^{lie}$ is a Lie algebra and $L^{ann}$ is a module over the Lie algebra 
$L^{lie}$ under the Lie module action:
\begin{equation}\label{eq1.5}
\overline{x}\cdot a:=-\langle a, x\rangle \qquad x\in\mathcal{L}, \qquad 
\end{equation}
where $x\in\mathcal{L}$, $\overline{x}:=x+ L^{ann}\in L^{lie}$ and $a\in L^{ann}$.

\medskip
Following \cite{LT}, we now define a module over a Leibniz algebra.

\begin{definition}\label{def1.5} Let $(L, \, \langle \, , \, \rangle )$ be a Leibniz algebra over a field $\mathbf{k}$. A vector space $V$ over the field $\mathbf{k}$ is called a {\bf module over the Leibniz algebra} (or a {\bf $L$-module}) if there are two linear maps 
$$ f: L\to End (L) \quad\mbox{and}\quad  g: L\to End (L)$$
such that
\begin{eqnarray}
\label{eq1.28} f(\langle x, y\rangle )&=& [f(x), f(y)],\\
\label{eq1.29} g(\langle x, y\rangle )&=& [g(x), f(y)],\\
\label{eq1.30} g(x)g(y)&=& g(x)f(y)
\end{eqnarray}
\end{definition}
for all $x$, $y\in L$. A $L$-module $V$ is also denoted by $(V, f, g)$, and $f$ and $g$ are called the {\bf right linear map} and the {\bf left linear map} of $V$, respectively.

\medskip
If $L$ is a Leibniz algebra, then $(L, -r, \ell)$ is a module over the Leibniz algebra $L$, where 
$$
r: L\to End(L) \quad\mbox{and}\quad \ell : L\to End(L),
$$
are the two linear maps defined by
\begin{equation}\label{eq1.18}
r_x(a):=\langle a, x\rangle ,\quad \ell _x(a):=\langle x, a\rangle \qquad\mbox{for $x$, $a\in\mathcal{L}$.}
\end{equation}
$r$ and $\ell$ are called the {\bf right multiplication} and the {\bf left multiplication} of $L$, respectively. This $L$-module $(L, -r, \ell)$ is called the {\bf adjoint module} over the Leibniz algebra $L$.

\medskip
Let $(V, f, g)$ be a module over a Leibniz algebra $L$. A subspace $U$ of $V$ is called a 
{\bf submodule} of $V$ if 
\begin{equation}\label{eq1.49}
f(x)(U)\subseteq U \quad\mbox{and}\quad g(x)(U)\subseteq U  \quad\mbox{for all $x\in L$.}
\end{equation}
The subspace
\begin{equation}\label{eq1.50}
V^{ann}:=\displaystyle\sum _{y\in L, \, v\in V} \mathbf{k}\Big(g(y)-f(y)\Big)(v)
\end{equation}
is called the {\bf annihilator} of $V$.

\medskip
The following proposition is clearly true.

\begin{proposition}\label{pr1.1} Let $(L, \, \langle \, , \, \rangle )$ be a Leibniz algebra over a field $\mathbf{k}$, and let $(V, f, g)$ be a module over the Leibniz algebra $L$.
\begin{description}
\item[(i)] The annihilator $V^{ann}$ is a Leibniz submodule of $V$.
\item[(ii)] The annihilator $V^{ann}$ is a module over the Lie algebra $\displaystyle\frac{L}{L^{ann}}$ under the following module action:
\begin{equation}\label{eq1.51}
\overline{a} \cdot v:= f(a)(v),
\end{equation}
where $\overline{a}:=a+L^{ann}\in \displaystyle\frac{L}{L^{ann}}$ and $v\in V^{ann}$.
\end{description}
\end{proposition}

\hfill\raisebox{1mm}{\framebox[2mm]{}}

\bigskip
We now introduce $3$-irreducible modules over Leibniz algebras.

\begin{definition}\label{def1.7} Let $(V, f, g)$ be a module over the Leibniz algebra $L$.
$(V, f, g)$ is said to be {\bf $3$-irreducible} if $V^{ann}\ne 0$ and $V$ has no submodules which are not equal to $0$, $V^{ann}$ and $V$.
\end{definition}

$3$-irreducible modules are  natural building blocks among the modules over Leibniz algebras. Finding the classification of $3$-irreducible modules over a Leibniz algebra is an interesting question, and will establish a representation theory which lifts the representation theory of 
semisimple Lie algebras.

\section{Root Space Decomposition}

A Leibniz algebra $(L, \, \langle \, , \, \rangle )$ is said to be {\bf solvable} if there exists a positive integer $n$ such that $\mathcal{D}^n L=0$, where $\mathcal{D}^n L$ is defined inductively as follows:
\begin{equation}\label{eq2.1}
\mathcal{D}^1 L:= L, \quad \mathcal{D}^{n+1} L=
\langle\mathcal{D}^n L, \,\mathcal{D}^n L\rangle\quad\forall n\ge 1.
\end{equation}
$L$ is said to be {\bf right nilpotent} if there exists a positive integer $n$ such that $\mathcal{C}_r^n L=0$, where $\mathcal{C}_r^n L$ is defined inductively by
\begin{equation}\label{eq3.1}
\mathcal{C}_r^1 L :=L, \quad \mathcal{C}_r^{n+1} L :=\langle \mathcal{C}_r^n L , L \rangle .
\end{equation}

\medskip
The Extended Lie's Theorem, which was announced in \cite{Liu1}, plays an important role in our study of root space decomposition.

\begin{proposition}\label{pr2.5} {\bf (The Extended Lie's Theorem)} Let $L$ be a finite dimensional solvable Leibniz algebra over an algebraically closed field of characteristic $0$. If 
$(V, f , g )$ is a finite dimensional nonzero $L$-module, then there exist two
linear functionals $\phi$, $\psi: L\to \mathbf{k}$ and a nonzero vector $v$ in $V$ such that
$$ f(z) (v)=\phi (z)v, \quad g(z)(v)=\psi (z)v \quad\mbox{for $z\in L$}. $$
Moreover, either $(\phi, \psi )=(\phi, \phi )$ or $(\phi, \psi )=(\phi, 0)$.
\end{proposition}

\hfill\raisebox{1mm}{\framebox[2mm]{}}

Using the Extended Lie's Theorem, we have

\begin{proposition}\label{pr4.5} Let $L$ be a right nilpotent Leibniz algebra over an algebraically closed field $\mathbf{k}$ of characteristic $0$, and let $(V, f, g)$ be a finite dimensional nonzero $L$-module.
\begin{description}
\item[(i)] $f(L)$ is a nilpotent Lie subalgebra of the general linear Lie algebra $g\ell (V)$ and
$V=\displaystyle\bigoplus _{\phi \in \overline{\Delta _f}}V^\phi _{f(L)}$, where  $\overline{\Delta _f}$ is the set of all weights of $V$ for $f(L)$, and $V^\phi _{f(L)}$ is the weight space of $V$ corresponding to the weight $\phi$; that is,
\begin{equation}\label{eq4.21}
0\ne V^\phi _{f(L)}:=\left\{\, v\in V \,\left |\,\begin{array}{c}\mbox{for every $x\in L$ there exists a positive}\\ \mbox{integer $n$ depending on $x$ and $v$}\\
\mbox{such that $\Big(f(x)-\phi (x)I\Big)^n(v)=0$}\end{array} \right. \, \right\}.
\end{equation}
\item[(ii)] For $\phi\in \overline{\Delta _f}$, $V^\phi _{f(L)}$ is a module over the right nilpotent Leibniz algebra $L$. Moreover, if $\phi\in \overline{\Delta _f}\setminus \{\, 0\,\}$, then $V^{ann}\bigcap V^\phi _{f(L)}=\Big(V^\phi _{f(L)}\Big)^{ann}$.
\item[(iii)] For $\phi \in \overline{\Delta _f}$, $V^\phi _{f(L)}$ has a basis:
\begin{equation}\label{eq4.37}
\underbrace{v_{11}, \cdots , v_{1n_1}}_{\mbox{a basis of $\Big(V^\phi _{f(L)}\Big)^{ann}$}}, v_1,  
v_{21}, \cdots , v_{2n_2}, v_2, \cdots 
\end{equation}
such that
\begin{equation}\label{eq4.38}
V^\phi _{f(L)}=\displaystyle\bigoplus _{i=1}^{p-1}
\left(\Big(\displaystyle\bigoplus _{j=1}^{n_i}\mathbf{k}v_{ij}\Big)\oplus \mathbf{k}v_i\right)
\oplus \Big(\displaystyle\bigoplus _{j=1}^{n_p}\mathbf{k}v_{pj}\Big)
\end{equation}
or
\begin{equation}\label{eq4.39}
V^\phi _{f(L)}=\displaystyle\bigoplus _{i=1}^p
\left(\Big(\displaystyle\bigoplus _{j=1}^{n_i}\mathbf{k}v_{ij}\Big)\oplus \mathbf{k}v_i\right)
\end{equation}
and the following hold:
\begin{equation}\label{eq4.40}
g(x)(v_{ih})\,\equiv\, 0\,mod\,\displaystyle\bigoplus _{s=1}^{i-1}
\left(\Big(\displaystyle\bigoplus _{j=1}^{n_s}\mathbf{k}v_{sj}\Big)\oplus \mathbf{k}v_s\right),
\end{equation}
\begin{equation}\label{eq4.41}
f(x)(v_i)\,\equiv\, \phi (x)v_i\,mod\,\left(\displaystyle\bigoplus _{s=1}^{i-1}
\Big(\Big(\displaystyle\bigoplus _{j=1}^{n_s}\mathbf{k}v_{sj}\Big)\oplus \mathbf{k}v_s\Big)\oplus \Big(\displaystyle\bigoplus _{j=1}^{n_i}\mathbf{k}v_{ij}\Big)\right),
\end{equation}
\begin{equation}\label{eq4.42}
g(x)(v_i)\,\equiv\, \phi (x)v_i\,mod\,\left(\displaystyle\bigoplus _{s=1}^{i-1}
\Big(\Big(\displaystyle\bigoplus _{j=1}^{n_s}\mathbf{k}v_{sj}\Big)\oplus \mathbf{k}v_s\Big)\oplus \Big(\displaystyle\bigoplus _{j=1}^{n_i}\mathbf{k}v_{ij}\Big)\right),
\end{equation}
where $x\in L$, $n_1$, $\cdots$, $n_p$ are nonnegative integers, 
$1\le i\le p$, $1\le h\le n_i$, and $\displaystyle\bigoplus _{j=1}^0\,(\star) :=0$.
\end{description}
\end{proposition}

The following definition gives the notion of a Cartan subalgebra of a Leibniz algebra, which has appeared in \cite{AAO}.

\begin{definition}\label{def2.1} Let $(L, \langle \, , \, \rangle )$ be a Leibniz algebra. A Leibniz subalgebra $H$ of $L$ is called a {\bf Cartan subalgebra} if the following two conditions are satisfied:
\begin{description}
\item[(i)] $H$ is right nilpotent;
\item[(ii)] $H$ coincides with the {\bf right normalizator} $N_r (H)$ of $H$, where  the right normalizator $N_r (H)$ of $H$ is defined by 
$N_r (H):=\{\, x\in L \,|\, \langle x, H\rangle \subseteq H \,\}$.
\end{description}
\end{definition}

Note that the set $\{\, x\in L \,|\, \langle x, H\rangle \subseteq H \,\}$ is called the left normalizator of $H$ in \cite{AAO}.

\medskip
Let $H$ be a Cartan subalgebra of a finite dimensional Leibniz algebra $L$. Using Proposition~\ref{pr4.5}, we have the following weight space decomposition of the adjoint module 
$(L, -r, \ell)$:
\begin{equation}\label{eq5.1}
L=L^0_{-r(H)}\bigoplus 
\left(\displaystyle\bigoplus _{\phi \in \Delta }L^\phi _{-r(H)}\right),
\end{equation}
where $\Delta$ is the set of all nonzero weights of $L$ for the linear nilpotent Lie algebra $-r(H)$. By 
(\ref{eq4.21}), we have
\begin{equation}\label{eq5.2}
L^0_{-r(H)}=\left\{\, y\in L \,\left |\,\begin{array}{c}\mbox{for every $h\in H$ there exists a positive}\\ \mbox{integer $n$ depending on $h$ and $y$}\\
\mbox{such that $(r_h)^n(y)=0$}\end{array} \right. \, \right\}.
\end{equation}
Using Proposition 2.3 in \cite{AAO} and Proposition 3.2 in \cite{O}, we know that 
$H=L^0_{-r(H)}$ and $\displaystyle\frac{H+L^{ann}}{L^{ann}}$ is a Cartan subalgebra of the Lie algebra $\displaystyle\frac{L}{L^{ann}}$. Thus, we can decompose $L$ as 
\begin{equation}\label{eq5.4}
L=H\bigoplus \left(\displaystyle\bigoplus _{\phi \in \Delta }L^\phi _{-r(H)}\right).
\end{equation}
The decomposition (\ref{eq5.4}) is called the {\bf root space decomposition} of $L$ with respect to the Cartan subalgebra $H$, the elements of $\Delta$ are called the {\bf roots} of $L$ with respect to the Cartan subalgebra $H$, and $L^\phi _{-r(H)}$ with $\phi \in\Delta$ is called the {\bf root space} of $L$ corresponding to the root $\phi$.

\medskip
We now introduce the notions of even roots and odd roots.

\begin{definition}\label{def5.2} Let $H$ be a Cartan subalgebra of a finite dimensional Leibniz algebra $L$, and let $L^\phi _{-r(H)}$ be a root space of $L$ corresponding to the root $\phi$ of $L$.
\begin{description}
\item[(i)] $L^\phi _{-r(H)}$ is called an {\bf even root space} and $\phi$ is called an 
{\bf even root} if $L^\phi _{-r(H)}\ne \Big(L^\phi _{-r(H)}\Big)^{ann}$.
\item[(ii)] $L^\phi _{-r(H)}$ is called an {\bf odd root space} and $\phi$ is called an {\bf odd root} if $L^\phi _{-r(H)}= \Big(L^\phi _{-r(H)}\Big)^{ann}$.
\end{description}
\end{definition}

Using the root space decomposition (\ref{eq5.4}), we have
\begin{equation}\label{eq5.16}
L=H\oplus \left(\displaystyle\bigoplus _{\phi \in \Delta _0}L^\phi _{-r(H)}\right)
\oplus \left(\displaystyle\bigoplus _{\phi \in \Delta _1}L^\phi _{-r(H)}\right),
\end{equation}
where
$$
\Delta _0:=\left\{\, \phi \,\left| \, \mbox{$\phi$ is a root of $L$ and 
$L^\phi _{-r(H)}\ne \Big(L^\phi _{-r(H)}\Big)^{ann}$}\right. \,\right\}
$$
is the set of all even roots of $L$ and 
$$
\Delta _1:=\left\{\, \phi \,\left| \, \mbox{$\phi$ is a root of $L$ and 
$L^\phi _{-r(H)}= \Big(L^\phi _{-r(H)}\Big)^{ann}$}\right. \,\right\}
$$
is the set of all odd roots of $L$.

\medskip
Using Proposition~\ref{pr4.5} again, we get the following proposition, which is useful in the study of simple Leibniz algebras.

\begin{proposition}\label{pr5.3} Let $L$ be a finite dimensional Leibniz algebra such that  the Lie-factor $L^{lie}=\displaystyle\frac{L}{L^{ann}}$ is a semisimple Lie algebra, let $H$ be a Cartan subalgebra of $L$, and let $L^\phi _{-r(H)}$ be a  root space of $L$ corresponding to a root $\phi$ of $L$.
\begin{description}
\item[(i)] $H$ is Abelian; that is $\langle H, H\rangle =0$. Moreover, we have 
\begin{equation}
L^{ann}=\Big(H\bigcap L^{ann}\Big)\bigoplus\left(\displaystyle\bigoplus _{\phi\in \Delta}
\Big(L^{\phi}_{-r(H)}\Big)^{ann}\right).
\end{equation}
\item[(ii)] If $\phi$ is an odd root, then there exist $n_{\phi}\in\mathcal{Z}_{\ge 1}$ and $v_1^{\phi}$, $\cdots$, $v_{n_{\phi}}^{\phi}\in L^\phi _{-r(H)}$ such that
\begin{equation}\label{eq5.17}
L^\phi _{-r(H)}=\mathbf{k}v_1^{\phi}\oplus \cdots \oplus v_{n_{\phi}}^{\phi}
\end{equation}
and the vectors $v_1^{\phi}$, $\cdots$, $v_{n_{\phi}}^{\phi}$ satisfy
\begin{equation}\label{eq5.18}
-r_h(v_i^{\phi})=\phi (h) v_i^{\phi}, \quad \ell _h(v_i^{\phi})=0 \quad\mbox{for $h\in H$ and $1\le i\le n_{\phi}$}.
\end{equation}
\item[(iii)] If $\phi$ is an even root, then there exist $n_{\phi}\in\mathcal{Z}_{\ge 0}$ and $v_1^{\phi}$, $\cdots$, $v_{n_{\phi}}^{\phi}$, $v^{\phi}\in L^\phi _{-r(H)}$ with $v^{\phi}\ne 0$ such that
\begin{equation}\label{eq5.19}
L^\phi _{-r(H)}=\mathbf{k}v_1^{\phi}\oplus \cdots \oplus \mathbf{k}v_{n_{\phi}}^{\phi}\oplus 
\mathbf{k}v^{\phi},
\end{equation}
the vectors $v_1^{\phi}$, $\cdots$, $v_{n_{\phi}}^{\phi}$ satisfy (\ref{eq5.18}) and the nonzero vector $v^{\phi}$ satisfies 
\begin{equation}\label{eq5.20}
-r_h(v^{\phi})=\phi (h) v^{\phi} \quad\mbox{for $h\in H$}
\end{equation}
and
\begin{equation}\label{eq5.21}
\ell _h(v^{\phi})\,\equiv\, \phi (h)(v^{\phi})\,mod\,\displaystyle\bigoplus _{i=1}^{n_{\phi}}
\mathbf{k}v_i^{\phi}  \quad\mbox{for $h\in H$},
\end{equation}
where $\mathbf{k}v_1^{\phi}\oplus \cdots \oplus \mathbf{k}v_{n_{\phi}}^{\phi}:=0$ for $n_{\phi}=0$.
\end{description}
\end{proposition}

\section{Simple Leibniz Algebras with Lie Factor $s\ell _2$}

Based on Proposition 1.2 in \cite{hl1}, we now define the building blocks in the context of Leibniz algebras by adding one more condition in the definition of simple Leibniz algebras introduced in \cite{D}.

\begin{definition}\label{def1.4} A Leibniz algebra $L$ is said to be {\bf simple} if 
$\langle L^{ann}, L \rangle \ne 0$ and $L$ has no ideals which are not equal to $\{0\}$, $L^{ann}$ and $L$.
\end{definition}

A Leibniz algebra $L$ is said to be {\bf standard} if there exists a Lie subalgebra $\bar{L}$ of $L$ such that $L=\bar{L}\oplus L^{ann}$ as a direct sum of vector spaces.

\medskip
We now construct two kinds of simple Leibniz algebras $s\ell _{2,2}^{\alpha, \beta}$ and 
$s\ell _{2,n\ge 4}^{\alpha, \beta}$,  where $\alpha$ and $\beta$ are two fixed elements of the field $\mathbf{k}$, and $n$ is an even integer with $n\ge 4$.

\medskip
\begin{quote}
$\bullet$ $s\ell _{2,2}^{\alpha, \beta}$ is a six-dimensional vector space:
$$
s\ell _{2,2}^{\alpha, \beta}=\mathbf{k}f\oplus\mathbf{k}h\oplus\mathbf{k}e\oplus\mathbf{k}v_0
\oplus\mathbf{k}v_1\oplus\mathbf{k}v_2 ,
$$ 
and the angle bracket on $s\ell _{2,2}^{\alpha, \beta}$ is defined by
$$
\langle h, e\rangle =2e+2\alpha v_0, \quad \langle e, h\rangle =-2e, 
\quad \langle h, f\rangle =-2f+\beta v_2,
$$
$$
\langle f, h\rangle =2f, \quad \langle e, f\rangle =h+\alpha v_1, 
\quad \langle f, e\rangle =-h-\beta v_1,
$$
$$
\langle v_k, h\rangle =-(2-2k)v_k, \quad \langle v_k, f\rangle =-v_{k+1}, \quad 
\langle v_k, e\rangle =-k(3-k)v_{k-1}
$$
and the other angle brackets of two elements in  the basis 
$$\{\, f, h, e, v_0, v_1, v_2 \,\}$$ 
are zero, where $k=0$, $1$, $2$ and $v_j:=0$ for $j\not \in \{\, 0, 1, 2\,\}$.
\end{quote}

\medskip
\begin{quote}
$\bullet$ $s\ell _{2,n\ge 4}^{\alpha, \beta}$ is a $(n+4)$-dimensional vector space:
$$
s\ell _{2,n\ge 4}^{\alpha, \beta}=\mathbf{k}f\oplus\mathbf{k}h\oplus\mathbf{k}e\oplus\mathbf{k}v_0
\oplus\mathbf{k}v_1\oplus\cdots\oplus\mathbf{k}v_n
$$ 
and the angle bracket on $s\ell _{2,n\ge 4}^{\alpha, \beta}$ is defined by:
$$
\langle e, e\rangle =\alpha v_{-2+\frac n2}, \quad \langle f, f\rangle =\beta v_{2+\frac n2}, \quad \langle h, e\rangle =2e-\alpha v_{-1+\frac n2}, 
$$
$$
\langle e, h\rangle =-2e, \quad \langle h, f\rangle =-2f+
\left(\frac{-4\alpha}{n^2+2n-8}-2\beta \right)v_{1+\frac n2}, \quad \langle f, h\rangle =2f,
$$
$$
\langle e, f\rangle =h, \quad 
\langle f, e\rangle =-h+\left(\frac{4\alpha}{n^2+2n-8}+\frac{n(n+2)}{4}\beta \right)v_{\frac n2}, 
$$
$$
\langle v_k, h\rangle =-(n-2k)v_k, \quad \langle v_k, f\rangle =-v_{k+1}, \quad 
\langle v_k, e\rangle =-k(n-k+1)v_{k-1},
$$
and the other angle brackets of two elements in  the basis 
$$\{\, f, h, e, v_0, v_1, \cdots , v_n \,\}$$ 
are zero, where $k=0$, $1$, $\cdots$, $n$ and $v_j:=0$ for $j\not \in \{\, 0, 1, \cdots , n\,\}$.
\end{quote}

\medskip
It is clear that both $s\ell _{2,2}^{\alpha, \beta}$ and 
$s\ell _{2,n\ge 4}^{\alpha, \beta}$ are simple Leibniz algebras and their annihilators are given by
$$
\Big(s\ell _{2,2}^{\alpha, \beta}\Big)^{ann}=\mathbf{k}v_0\oplus\mathbf{k}v_1\oplus\mathbf{k}v_2 
$$
and
$$
\Big(s\ell _{2,n\ge 4}^{\alpha, \beta}\Big)^{ann}
=\mathbf{k}v_0\oplus\mathbf{k}v_1\oplus\cdots\oplus\mathbf{k}v_n .
$$

\medskip
Although there exists an odd-dimensional non-standard simple Leibniz algebra with Lie factor 
$s\ell _2$ in characteristic $p(>3)$ (\cite{DA}), any odd-dimensional simple Leibniz algebra with Lie factor $s\ell _2$ in characteristic $0$ is always standard by the following classification of finite dimensional simple Leibniz algebra with Lie factor $s\ell _2$.

\begin{proposition} Let $L$ be a finite dimensional simple Leibniz algebra over an algebraically closed field $\mathbf{k}$ of characteristic $0$. If the Lie factor $\displaystyle\frac{L}{L^{ann}}$ is the three-dimensional simple Lie algebra $s\ell _2$, then
\begin{description}
\item[(i)] $dim\, L$ is odd $\Longrightarrow$ $dim\, L\ge 5$ and $L$ is standard;
\item[(ii)] $dim\, L$ is even $\Longrightarrow$ $dim\, L\ge 6$ and there exist $\alpha$, $\beta\in\mathbf{k}$ such that
$$
\mbox{$L$ is isomorphic to}\left\{\begin{array}{l} \mbox{$s\ell _{2,2}^{\alpha, \beta}$ for 
$dim\, L=6$}\\ \\ \mbox{$s\ell _{2,n\ge 4}^{\alpha, \beta}$ for 
$dim\, L > 6$}\end{array}\right. ,
$$
where $n=dim\, L -4$.
\end{description}
\end{proposition}

\hfill\raisebox{1mm}{\framebox[2mm]{}}

\end{document}